\documentclass[10pt]{article}
\usepackage{amssymb}
\usepackage{amsthm}
\usepackage{graphicx}
\usepackage{graphics}
\usepackage{amsthm}

\input pictex.tex
\newcommand{\esp}{\hspace{0.1cm}}

\newcommand{\efe}{\mathrm{F}}

\theoremstyle{definition}

\newtheorem{thm}{Theorem}[section]

\newtheorem{prop}[thm]{Proposition}

\newtheorem{lem}[thm]{Lemma}

\newtheorem{rem}[thm]{Remark}

\newtheorem{ex}[thm]{Example}

\newtheorem{qs}[thm]{Question}

\input epsf

\begin{document}

\date{}

\author{Andr\'es Navas \hspace{0.25cm} \& \hspace{0.25cm} Crist\'obal Rivas}

\title{Describing all bi-orderings on Thompson's group F}
\maketitle

\vspace{-0.25cm}

\noindent{\bf Abstract:} We describe all possible ways of bi-ordering Thompson's group F: 
its space of bi-orderings is made up of eight isolated points and four canonical copies 
of the Cantor set. 

\vspace{0.2cm}

\noindent{\bf Keywords:} Group orderings, Conrad property, Thompson's group F.

\vspace{0.2cm}

\noindent{\bf AMS Mathematics Subject Classification:} 06F15, 20F60, 22F50.

\vspace{0.65cm}

\noindent{\Large {\bf Introduction}}

\vspace{0.4cm}

In recent years, the well developed theory of orderable groups has re-emerged, 
mainly due to its connexions with many different branches of mathematics. One of the 
aspects which has been emphasized is that, in general, orderable groups actually admit 
many invariant total order relations. This makes natural the problem of searching for an 
ordering satisfying a nice property implying a relevant algebraic (or dynamical) property 
of the underlying group. This issue has been successfully exploited for instance by 
Witte-Morris in his beautiful proof of the local indicability for left-orderable 
amenable groups \cite{witte}. The reader is referred to \cite{ordering} for other 
applications of this approach.

A closely related problem concerns the description of {\em all} (invariant) orderings on particular 
classes of groups. In this direction, Tararin's concise classification of groups admitting only finitely 
many left-orderings corresponds to a relevant piece of the theory \cite{koko}. Another significant 
(and easier) result is the description of all possible orderings on torsion-free finite rank Abelian 
groups \cite{robbin,sikora,teh}.

Although the description of {\em all} orderings seems to be out of reach for general orderable 
groups, one may address the weaker question of the description of the corresponding space of 
orderings from a topological viewpoint (recall that the space of orderings on 
any space corresponds to the projective limit of the orders on finite sets, and 
hence carries the structure of a compact topological space). For instance, ruling 
out the existence of {\em isolated points} in this space (that is, orderings which are completely 
determined by finitely many inequalities) appears to be a fundamental question. This has been done 
for instance for the spaces of left-orderings of finitely generated torsion-free nilpotent groups 
which are not rank-1 Abelian \cite{ordering,sikora}. For the free group $F_n$ (where $n \geq 2$), 
it is known that there is no isolated point in the corresponding space of left-orderings 
\cite{MC,ordering,smith}. The similar question for the space of bi-orderings on $F_n$ 
remains open, and though it is not treated here, it inspires much of this work.  

\vspace{0.15cm}

In this article, we focus on a remarkable bi-orderable group, namely Thompson's group 
$\mathrm{F}$, and we provide a complete description of all its possible bi-orderings. 
Recall that $\mathrm{F}$ is the group of orientation-preserving piecewise-linear 
homeomorphisms $f$ of the interval $[0,1]$ such that:

\vspace{0.08cm}

\noindent -- the derivative of $f$ on each linearity interval is an integer power of $2$,

\vspace{0.08cm}

\noindent -- $f$ induces a bijection of the set of dyadic rational numbers in $[0,1]$.

\vspace{0.08cm}

\noindent For each non-trivial $f \!\in\! \mathrm{F}$ we will denote by $x^{-}_{f}$ (resp. 
$x^{+}_f$) the leftmost point $x^-$ (resp. the rightmost point $x^+$) for which $f'_+ (x^{-}) \neq 1$ 
(resp. $f'_{-} (x^{+}) \neq 1$), where $f'_+$ and $f'_{-}$ stand for the corresponding lateral 
derivatives. One can then immediately visualize four different bi-orderings on (each subgroup of) 
$\mathrm{F}$, namely:\\

\vspace{0.1cm}

\noindent -- the bi-ordering $\preceq_{x^{-}}^{+}$ for which \esp $f \succ id $ 
\esp \esp if and only if \esp $f'_+ (x^{-}_f) > 1$,\\ 

\vspace{0.1cm}

\noindent -- the bi-ordering $\preceq_{x^{-}}^{-}$ for which \esp $f \succ id $ 
\esp \esp if and only if \esp $f'_+ (x^{-}_f) < 1$,\\ 

\vspace{0.1cm}

\noindent -- the bi-ordering $\preceq_{x^{+}}^{+}$ for which \esp $f \succ id $ 
\esp \esp if and only if \esp $f'_{-} (x^{+}_f) < 1$,\\ 

\vspace{0.1cm}

\noindent -- the bi-ordering $\preceq_{x^{+}}^{-}$ for which \esp $f \succ id $ 
\esp \esp if and only if \esp $f'_{-} (x^{+}_f) > 1$.\\

\vspace{0.1cm}

\noindent Although $\mathrm{F}$ admits many more bi-orderings than these, 
the case of its derived subgroup $\mathrm{F}'$ is quite different.

\vspace{0.5cm}

\noindent{\bf Theorem [V. Dlab].} {\em The only bi-orderings on $\mathrm{F}'$ 
are $\preceq_{x^{-}}^{+}$, $\preceq_{x^{-}}^{-}$, $\preceq_{x^{+}}^{+}$ and 
$\preceq_{x^{+}}^{-}$.}

\vspace{0.5cm}

Dlab's arguments apply to many other (in general, non finitely generated) groups of 
piecewise-affine homeomorphisms of the line. Some of them appear to be non-Abelian, 
though having only two different bi-orderings (compare Remark \ref{cuatro}). We refer 
to the original reference \cite{dlab} for all of this (see also \cite{koko,Me1,Me2,ZM}). 
Here we provide a new proof using an argument which allows us to obtain the complete 
classification of all the bi-orderings on $\mathrm{F}$.  

\vspace{0.1cm}

Remark that there are also four other ``exotic'' bi-orderings on F, namely:\\

\vspace{0.1cm} 

\noindent -- the bi-ordering $\preceq_{0,x^{-}}^{+,-}$ for which $f \succ id$ if and only 
if either $x^{-}_f = 0$ and $f'_+(0) > 1$, or $x^{-}_f \neq 0$ and $f'_+ (x^{-}_f) < 1$,\\ 

\vspace{0.1cm}

\noindent -- the bi-ordering $\preceq_{0,x^{-}}^{-,+}$ for which $f \succ id$ if and only 
if either $x^{-}_f = 0$ and $f'_+(0) < 1$, or $x^{-}_f \neq 0$ and $f'_+ (x^{-}_f) > 1$,\\ 

\vspace{0.1cm}

\noindent -- the bi-ordering $\preceq_{1,x^{+}}^{+,-}$ for which $f \succ id $ if and only 
if either $x^{+}_f = 1$ and $f'_+(1) < 1$, or $x^{+}_f \neq 1$ and $f'_{-} (x^{+}_f) > 1$,\\ 

\vspace{0.1cm}

\noindent -- the bi-ordering $\preceq_{1,x^{+}}^{-,+}$ for which $f \succ id $ if and only 
if either $x^{+}_f = 1$ and $f'_+(1) > 1$, or $x^{+}_f \neq 1$ and $f'_{-} (x^{+}_f) < 1$.\\

\vspace{0.1cm}

\noindent Notice that, when restricted to $\mathrm{F}'$, the bi-ordering $\preceq_{0,x^{-}}^{+,-}$ 
(resp. $\preceq_{0,x^{-}}^{-,+}$, $\preceq_{1,x^{+}}^{+,-}$, and $\preceq_{1,x^{+}}^{-,+}$) 
coincides with $\preceq_{x^{-}}^{-}$ (resp. $\preceq_{x^{-}}^{+}$, $\preceq_{x^{+}}^{-}$, 
and $\preceq_{x^{+}}^{+}$). Let us denote the set of the previous eight bi-orderings on F 
by $\mathcal{BO}_{Isol} (\mathrm{F})$.

\vspace{0.15cm}

There is another natural procedure for creating bi-orderings on F. For this, recall the 
well-known (and easy to check) fact that $\mathrm{F}'$ coincides with the subgroup of 
$\mathrm{F}$ formed by the elements $f$ satisfying $f'_+(0) = f'_{-}(1) = 1$. Now  
let $\preceq_{\mathbb{Z}^2}$ be any bi-ordering on $\mathbb{Z}^2$, and let 
$\preceq_{\mathrm{F}'}$ be any bi-ordering on $\mathrm{F}'$. It readily follows from 
Dlab's theorem that 
$\preceq_{\efe'}$ is invariant under conjugacy by elements in $\efe$. Hence,  
one may define a bi-ordering $\preceq$ on F by declaring that 
$f \succ id$ if and only if either $f \notin \mathrm{F}'$ 
and $\big( \log_2 (f'_{+}(0)), \log_2 (f'_{-}(1))\big) 
\succ_{\mathbb{Z}^2} \big( 0,0 \big)$, or $f \in \mathrm{F}'$ 
and $f \succ_{\mathrm{F}'} id$.

All possible ways of ordering finite-rank Abelian groups have been described in \cite{robbin,sikora,teh}. 
In particular, when the rank is greater than 1, the corresponding spaces of bi-orderings are 
homeomorphic to the Cantor set. Since there are only four possibilities for the 
bi-ordering $\preceq_{\mathrm{F}'}$, the preceding procedure gives four natural copies 
(which we will coherently denote by $\Lambda_{x^{-}}^{+}$, $\Lambda_{x^{-}}^{-}$, 
$\Lambda_{x^{+}}^{+}$, and $\Lambda_{x^{+}}^{-}$) of the Cantor set in the space of bi-orderings 
of F. The main result of this work establishes that these bi-orderings, together with the special 
eight bi-orderings previously introduced, fill out the list of all possible bi-orderings on F.

\vspace{0.5cm}

\noindent{\bf Theorem.} {\em The space of bi-orderings of \esp $\mathrm{F}$ is the disjoint 
union of the finite set $\mathcal{BO}_{Isol} (\mathrm{F})$ (whose elements are isolated 
bi-orderings) and the copies of the Cantor set $\Lambda_{x^{-}}^{+}$, $\Lambda_{x^{-}}^{-}$, 
$\Lambda_{x^{+}}^{+}$, and $\Lambda_{x^{+}}^{-}$.}

\vspace{0.5cm}

The first ingredient of the proof of this result comes from the theory of Conradian orderings 
\cite{conrad}. Indeed, since $\mathrm{F}$ is finitely generated, every bi-ordering $\preceq$ 
on it admits a maximal proper convex subgroup $\mathrm{F}^{max}_{\preceq}$. More importantly, 
this subgroup may be detected as the kernel of a non-trivial, non-decreasing group homomorphism 
into $(\mathbb{R},+)$. Since $\mathrm{F}'$ is simple (see for instance \cite{CFP}) and 
non-Abelian, it must be contained in $\mathrm{F}^{max}_{\preceq}$. The case of coincidence 
is more or less transparent: the bi-ordering on $\mathrm{F}$ is contained in one of the four canonical 
copies of the Cantor set, and the corresponding bi-ordering on $\mathbb{Z}^2$ is of {\em irrational 
type} ({\em i.e.}, its positive elements are those which are in one of the two 
half-planes determined by a line of irrational slope passing through the origin). 
The case where $\mathrm{F}'$ is strictly contained in $\mathrm{F}^{max}_{\preceq}$ is more 
complicated. The bi-ordering may still be contained in one of the four canonical copies of the Cantor 
set, but the corresponding bi-ordering on $\mathbb{Z}^2$ must be of {\em rational type} ({\em e.g.}, 
a lexicographic ordering). However, it may also coincide with one of the eight special bi-orderings 
listed above. Distinguishing these two possibilities is the hardest part of the proof. For this, we 
strongly use the internal structure of $\mathrm{F}$, in particular the fact that the subgroup 
consisting of elements whose support is contained in a prescribed closed dyadic interval is 
isomorphic to $\mathrm{F}$ itself.

\vspace{0.5cm}

\noindent{\bf Acknowledgments.} The first author would like to thank Jos\'e Burillo for 
his explanations on the group of outer automorphisms of F, as well as Adam Clay and 
Dale Rolfsen for helpful discussions on orderable groups. This work was partially 
funded by the PBCT-Conicyt via the Research Network on Low Dimensional Dynamical 
Systems. The second author was also funded by the Conicyt PhD Fellowship 21080054.


\section{Some background}

\subsection{On group orderings}

\hspace{0.35cm} Throughout this work, the word {\em left-ordering} (resp. {\em bi-ordering}) will 
stand for a total order relation on a group which is invariant by left multiplication (resp. 
by left and right multiplication simultaneously). An element $f$ is said to be {\em positive} 
(resp. {\em negative}) with respect to some left-ordering $\preceq$ if $f \succ id$ (resp. 
$f \prec id$). The set of positive elements forms a semigroup $P_{\preceq}^+$, which is 
called the {\em positive cone} of $\preceq$, and the whole 
group equals the disjoint union of $P_{\preceq}^+$ together with 
$P_{\preceq}^{-} = \{f \! : f^{-1} \in P_{\preceq}^+ \}$ and $\{ id \}$. Conversely, 
given a subsemigroup $P^+$ of a group $\Gamma$ such that $\Gamma$ equals the disjoint 
union of $P^{+}$ together with $P^{-} \!= \{f \! : f^{-1} \in P^{+}\}$ and $\{id\}$, one 
may realize $P^{+}$ as the positive cone of a left-ordering $\preceq$: it suffices to declare 
that $f \succ g$ if and only if $g^{-1} f$ belongs to $P^{+}$. The resulting ordering will be 
bi-invariant if and only if $P^{+}$ is a normal subsemigroup, that is, if $gfg^{-1} \in P^{+}$ 
for all $f \in P^{+}$ and all $g \in \Gamma$.

Every left-ordering (resp. bi-ordering) $\preceq$ on a group $\Gamma$ comes together with an 
associated ({\em conjugate}) left-ordering (resp. bi-ordering) $\bar{\preceq}$ whose positive 
cone coincides with $P_{\preceq}^-$. Clearly, the map \esp $\preceq \esp \mapsto \bar{\preceq}$ 
\esp is an involution of the set of left-orderings (resp. bi-orderings).

\begin{ex} Clearly, there are only two bi-orderings on $\mathbb{Z}$. The case of $\mathbb{Z}^2$ is 
more interesting. According to \cite{robbin,sikora,teh}, there are two different types of bi-orderings 
on $\mathbb{Z}^2$. Bi-orderings of {\em irrational type} are completely determined by an irrational number 
$\lambda$: for such an order $\preceq_{\lambda}$ an element $(m,n)$ is positive if and only if \esp 
$\lambda m + n$ \esp is a positive real number. Bi-orderings of rational type are characterized by 
two data, namely a pair $(a,b) \in \mathbb{Q}^2$ up to multiplication by a positive real number, 
and the choice of one of the two possible bi-orderings on the subgroup \esp 
$\{(m,n) \!: am + bn = 0\} \sim \mathbb{Z}$. \esp Thus an element 
\esp $(m,n) \in \mathbb{Z}^2$ \esp is 
positive if and only if either \esp $am + bn$ \esp 
is a positive real number, or \esp $am+bn = 0$ \esp and 
$(m,n)$ is positive with respect to the chosen bi-ordering on the kernel line 
(isomorphic to $\mathbb{Z}$). The description of all bi-orderings on $\mathbb{Z}^n$ 
for bigger $n$ continues inductively. A good exercise is to obtain all of this 
by using Conrad's theorem from \S \ref{con}.
\label{zeta2}
\label{abel}
\end{ex}


\subsection{On spaces of orderings}

\hspace{0.35cm} Given a left-orderable group $\Gamma$ (of arbitrary cardinality), 
we denote by $\mathcal{LO} (\Gamma)$ the set of all left-orderings on $\Gamma$. 
This set has a natural topology: a basis of neighborhoods of $\preceq$ in 
$\mathcal{LO} (\Gamma)$ is the family of the sets $U_{g_0,\ldots,g_k}$ of all 
left-orderings $\preceq'$ on $\Gamma$ which coincide with $\preceq$ on $\{g_0,\ldots,g_k\}$, 
where $\{g_0,\ldots,g_k\}$ runs over all finite subsets of $\Gamma$. Endowed with this 
topology, $\mathcal{LO}(\Gamma)$ is totally disconnected, and by (an easy application 
of) the Tychonov Theorem, it is compact. The (perhaps empty) subspace $\mathcal{BO}(\Gamma)$ 
of bi-orderings on $\Gamma$ is closed inside $\mathcal{LO}(\Gamma)$, and hence is also compact. 

If $\Gamma$ is countable, then the above topology is metrizable: given an exhaustion 
$\Gamma_0 \subset \Gamma_1 \subset \ldots$ of $\Gamma$ by finite sets, for different 
$\preceq$ and $\preceq'$ one may define $dist (\preceq,\preceq') = 1 / 2^n$, where $n$ 
is the first integer such that $\preceq$ and $\preceq'$ do not coincide on $\Gamma_n$. 
If $\Gamma$ is finitely generated, one may take $\Gamma_n$ as being the ball of radius 
$n$ with respect to some fixed finite system of generators. (The metrics arising from 
two different finite systems of generators are H\"older equivalent.) 

By definition, an isolated point $\preceq$ in $\mathcal{LO}(\Gamma)$ corresponds to an 
ordering for which there exist $g_0,\ldots,g_k$ in $\Gamma$ such that $U_{g_0,\ldots,g_k}$ 
reduces to $\{ \preceq \}$. This is the case for example if $g_1,\ldots,g_k$ generate the positive 
cone of $\preceq$ as a semigroup and $g_0 = id$: see \cite[Proposition 1.8]{ordering}. 
Analogously, $\preceq$ is an isolated point of $\mathcal{BO}(\Gamma)$ if 
$U_{g_0,\ldots,g_n} \cap \mathcal{BO}(\Gamma)$ reduces to $\{ \preceq \}$ for some 
$g_0,\ldots,g_k$ in $\Gamma$. According to the (obvious) proposition below, this happens 
for instance if $g_1,\ldots,g_k$ generate the positive cone of $\preceq$ as a 
{\em normal} semigroup and $g_0 = id$ (recall that a subset $S$ of a normal 
subsemigroup $P$ of a group $\Gamma$ generates $P$ as a normal semigroup if 
$P$ coincides with the smallest normal subsemigroup $\langle S \rangle_N^+$ 
of $\Gamma$ containing $S$): see Questions \ref{prima} and \ref{secua} on this.

\vspace{0.1cm}

\begin{prop} {\em Suppose that the positive cone of a bi-ordering $\preceq$ on a 
group $\Gamma$ is generated as a normal semigroup by elements $g_1,\ldots,g_k$. Then 
$\preceq$ is the unique bi-ordering on $\Gamma$ for which all of these elements are positive.}
\label{lema-obvio}
\end{prop}

\vspace{0.1cm}

As has been remarked by many people (see for instance \cite{ordering}), the group of 
automorphisms $Aut (\Gamma)$ of a left-orderable group $\Gamma$ acts by homeomorphisms 
of $\mathcal{LO} (\Gamma)$: given $\gamma \in Aut(\Gamma)$ and $\preceq$ in 
$\mathcal{LO}(\Gamma)$, the image of $\preceq$ by $\gamma$ is the left-ordering $\preceq_{\gamma}$  
whose positive cone is the preimage under $\gamma$ of the positive cone of $\preceq$. If $\Gamma$ 
is bi-orderable, then this action restricted to $\mathcal{BO} (\Gamma)$ factors through the 
group of outer automorphisms $Out (\Gamma)$. 

The dynamical properties of the preceding action for general bi-orderable groups seem interesting. 
For instance, the action of $GL(2,\mathbb{Z})$ on $\mathcal{BO} (\mathbb{Z}^2)$ is transitive on 
the set of bi-orderings of rational type, while the set of bi-orderings of irrational type 
decomposes into uncountably many orbits ({\em c.f.,} Example \ref{zeta2}). 

In a similar direction, the action of $Out (F_n)$ could be useful for understanding 
$\mathcal{BO} (\mathrm{F}_n)$. Nevertheless, in the case of Thompson's group $\mathrm{F}$, 
the action of $Out (\mathrm{F})$ on $\mathcal{BO} (\efe)$ is almost trivial. Indeed, 
according to \cite{ihes}, the group $Out (\mathrm{F})$ contains an index-two subgroup 
$Out_+ (\mathrm{F})$ whose elements are (equivalence classes of) conjugacies by certain 
orientation preserving homeomorphisms of the interval $[0,1]$. Although these homeomorphisms 
are dyadically piecewise-affine on $]0,1[$, the points of discontinuity of their derivatives 
may accumulate at $0$ and/or $1$, but in some ``periodically coherent'' way. It turns 
out that the conjugacies by these homeomorphisms preserve the derivatives of 
non-trivial elements $f \!\in\! \mathrm{F}$ at the points $x_f^-$ and $x_f^+$: this 
is obvious when these points are different from $0$ and $1$, and in the other case this 
follows from the explicit description of $Out (\mathrm{F})$ given in \cite{ihes}. According 
to our main theorem, this implies that the action of $Out_+ (\mathrm{F})$ on 
$\mathcal{BO}(\mathrm{F})$ is trivial. 

The set \esp $Out (\mathrm{F}) \setminus Out_+ (\mathrm{F})$ \esp corresponds to the class of 
the order-two automorphism $\sigma$ induced by the conjugacy by the map $x \mapsto 1-x$. One can 
easily check that 
$$(\preceq_{x^-}^+)_{\sigma} = \preceq_{x^+}^{-}, \quad (\preceq_{x^-}^-)_{\sigma} = \preceq_{x^+}^{+}, 
\quad (\preceq_{0,x^-}^{+,-})_{\sigma} = \preceq_{1,x^+}^{-,+}, \quad \mbox{and} \quad  
(\preceq_{0,x^-}^{-,+})_{\sigma} = \preceq_{1,x^+}^{+,-}.$$ 
Moreover, $\sigma (\Lambda_{x^-}^{+}) = \Lambda_{x^+}^{-}$ and 
$\sigma (\Lambda_{x^-}^{-}) = \Lambda_{x^+}^{+}$, and the action 
on the bi-orderings of the $\mathbb{Z}^2$-fiber can be easily 
described. We leave the details to the reader. 

\begin{rem} As in the case of $\sigma$, the dynamics of the involution 
$\preceq \esp \mapsto \bar{\preceq}$ can be also easily described. However, in the 
case of $\efe$, this involution does not occur as the action of any group automorphism.
\end{rem}


\subsection{On Conradian orderings}
\label{con}

\hspace{0.35cm} Besides  
$\mathcal{BO}(\Gamma)$, for a left-orderable group $\Gamma$ there is another relevant (perhaps 
empty) closed subset of $\mathcal{LO}(\Gamma)$, namely the subset $\mathcal{CO}(\Gamma)$ formed 
by the left-orderings $\preceq$ such that $g^{-1} f g^2 \succ id$ for all positive elements 
$f,g$ (see for instance \cite{conrad,ordering}). A left-ordering satisfying this property is 
said to be a $\mathcal{C}$-ordering or a {\em Conradian} ordering, and a group admitting 
such a left-ordering is called {\em Conrad-orderable} or simply $\mathcal{C}$-orderable.  
Notice that every bi-invariant ordering is Conradian.

In \cite{conrad}, a structure theory for Conradian orderings is given. 
(An alternative dynamical approach appears in \cite{ordering,crossing}.) 
This is summarized in the theorem below. To state it properly, recall that a 
subgroup $\Gamma_0$ of a group $\Gamma$ endowed with a left-ordering $\preceq$ is said to 
be $\preceq${\em-convex} if every $g \in \Gamma$ satisfying $g_{1} \preceq g \preceq g_2$ 
for some $g_1,g_2$ in $\Gamma_0$ actually belongs to $\Gamma_0$. Equivalently, every 
$h \in \Gamma$ satisfying $id \preceq h \preceq g$ for some $g \in \Gamma_0$ is contained 
in $\Gamma_0$. Notice that given any two $\preceq$-convex subgroups of $\Gamma$, one of 
them is necessarily contained in the other. Consequently, the union and the intersection 
of groups in an arbitrary family of $\preceq$-convex subgroups is also $\preceq$-convex.

\vspace{0.5cm}

\noindent{\bf Theorem [P. Conrad].} {\em Let $\Gamma$ be a group endowed with a $\mathcal{C}$-ordering. 
Given $g \in \Gamma$, denote by $\Gamma_g$ (resp. $\Gamma^g$) the maximal (resp. minimal) convex subgroup 
which does not contain (which contains) $g$. Then $\Gamma_g$ is normal in $\Gamma^g$, and there exists an 
non-decreasing group homomorphism $\tau_{\preceq}^{g} \!: \Gamma \rightarrow (\mathbb{R},+)$ whose kernel 
coincides with $\Gamma_g$. This homomorphism is unique up to multiplication by a positive real number. 

Moreover, if \esp $\Gamma$ is finitely generated, then it contains a (unique) maximal proper 
$\preceq$-convex subgroup $\Gamma^{max} = \Gamma^{max}_{\preceq}$, which coincides with 
the kernel of a (unique up to multiplication by a positive real number) non-decreasing 
group homomorphism $\tau_{\preceq} \!: \Gamma \rightarrow (\mathbb{R},+)$.}

\vspace{0.5cm}

A direct consequence of this theorem is that Conrad-orderable groups are locally indicable, that 
is, their  non-trivial finitely generated subgroups admit non-trivial group homomorphisms into 
$(\mathbb{R},+)$. Actually, the converse is also true (see for instance \cite{ordering} and 
references therein). 

The study of the topological properties of $\mathcal{CO}(\Gamma)$ is much simpler than those of 
$\mathcal{BO}(\Gamma)$. Indeed, in most of the cases, $\mathcal{CO}(\Gamma)$ has no isolated 
point (and hence it is homeomorphic to the Cantor set if the group is countable). To 
show a result in this direction, we need to recall the {\em extension procedure} 
for creating group orderings. 

Let $\preceq$ be a left-ordering on a group $\Gamma$, let $\Gamma_0$ be a $\preceq$-convex subgroup 
of $\Gamma$, and let $\preceq_0$ be a left-ordering on $\Gamma_0$. The extension of $\preceq_0$ by $\preceq$ 
is the left-ordering $\preceq^*$ on $\Gamma$ obtained by ``changing'' $\preceq$ into $\preceq_0$ on 
$\Gamma_0$, and ``keeping it'' outside. More precisely, the positive cone of $\preceq^*$ is 
$P^+_{\preceq_0} \cup (P^+_{\preceq} \setminus \Gamma_0)$. One can easily check that $\Gamma_0$ 
remains $\preceq^*$-convex. Moreover, if $\preceq$ and $\preceq_0$ are Conradian, then the 
resulting $\preceq^*$ is also a $\mathcal{C}$-ordering. Unfortunately (or perhaps 
fortunately), the  bi-invariance 
of both $\preceq$ and $\preceq_0$ does not guarantee the bi-invariance of $\preceq^*$: 
to ensure this, we also need to assume that the positive cone of $\preceq_0$ is invariant 
under conjugacies by elements in $\Gamma$.   
Finally, it is not difficult to check that if $\Gamma_0$ is a $\preceq$-convex {\em normal} 
subgroup of $\Gamma$, then $\preceq$ induces a left-ordering on the quotient 
$\Gamma / \Gamma_0$, which is a bi-ordering if $\preceq$ is bi-invariant.

\vspace{0.05cm}

\begin{ex} To simplify, denote just by $\preceq$ the bi-ordering $\preceq_{x^+}^{+}$ on $\mathrm{F}$. 
For a non-trivial element $g \in \mathrm{F}$, the subgroups $\mathrm{F}_g$ and $\mathrm{F}^g$ 
coincide with $\{f \in \mathrm{F} \! : supp(f) \subset ]x_g^{-},1] \}$ and 
$\{f \in \mathrm{F} \! : supp(f) \subset [x_g^{-},1] \}$ respectively, where 
$supp(f) = \overline{\{x \! : f(x) \neq x \}}$ is the {\em support} of $f$. The quotient 
$\Gamma^g / \Gamma_g$ is order isomorphic to $\mathbb{Z}$ via the homomorphism 
$f \Gamma_g \mapsto \log_2 \big( f'_+ (x_{g}^{-}) \big)$. A curious $\mathcal{C}$-ordering $
\preceq'$ on $\mathrm{F}$ (which is not bi-invariant~!) is obtained as follows: take the 
extension $\preceq^*$ of the restriction of $\preceq$ to $\Gamma_g$ by the restriction 
of $\bar{\preceq}$ to $\Gamma^g$, and then extend $\preceq^*$ by $\preceq$. This 
left-ordering obeys the following rule: a non-trivial element $f \!\in\! \mathrm{F}$ 
is positive with respect to $\preceq'$ if and only if either $x_f^{-} \neq x_g^{-}$ 
and $f'_+ (x_f^{-}) > 1$, or $x_f^{-} = x_g^{-}$ and $f'_+ (x_f^{-}) < 1$.
\end{ex}

\begin{ex} As the reader can easily check, the bi-ordering $\preceq_{0,x^{-}}^{+,-}$ 
appears as the extension by $\preceq_{x^{-}}^{+}$ of the restriction of its conjugate 
$\bar{\preceq}_{x^{-}}^{+}$ (which coincides with $\preceq_{x^{-}}^{-}$) to the maximal proper 
${\preceq_{x^{-}}^{+}}$-convex subgroup $\mathrm{F}^{max} = \{f \!\in\! \efe \!:\ f'_+(0) = 1\}$.  
The bi-orderings $\preceq_{0,x^{-}}^{-,+}$, $\preceq_{1,x^{+}}^{+,-}$, and $\preceq_{1,x^{+}}^{-,+}$  
may be obtained in the same way starting from $\preceq_{x^{-}}^{-}$, $\preceq_{x^{+}}^{+}$, and 
$\preceq_{x^{+}}^{-}$, respectively.
\end{ex}

\begin{rem}
In general, if $\Gamma$ is a finitely generated (non-trivial) group endowed with a 
bi-ordering $\preceq$, 
one can easily check that the ordering $\preceq^{*}$ obtained as the extension 
by $\preceq$ of $\bar{\preceq}$ 
restricted to $\Gamma^{max}_{\preceq}$ is bi-invariant. This bi-ordering 
(resp. its conjugate $\bar{\preceq}_*$) is always different from $\bar{\preceq}$ 
(resp. from $\preceq$), and it coincides with $\preceq$ (resp. with $\bar{\preceq}$) 
if and only if the only proper $\preceq$-convex subgroup is the trivial one; by Conrad's theorem, 
$\Gamma$ is necessarily Abelian in this case. We thus conclude that every non-Abelian 
finitely generated bi-orderable group admits at least four different bi-orderings. Moreover, 
(non-trivial) torsion-free Abelian groups having only two bi-orderings are those of rank 
one (in higher rank one may consider lexicographic type orderings).
\label{cuatro}
\end{rem}

\vspace{0.05cm}

\begin{prop} {\em If $\Gamma$ is a non-solvable Conrad-orderable group, 
then $\mathcal{CO}(\Gamma)$ contains no isolated point.}
\end{prop}

\noindent{\bf Proof.} Throughout the proof, fix a 
$\mathcal{C}$-ordering $\preceq$ on $\Gamma$. We will first show that if there 
are infinitely many subgroups of the form $\Gamma_g$, then $\preceq$ is not isolated inside 
$\mathcal{CO} (\Gamma)$. Indeed, given finitely many distinct elements $g_1,\ldots,g_k$ in 
$\Gamma$, consider the elements $f_{i,j}$ of the form $g_i^{-1} g_j$. We need to produce a  
$\mathcal{C}$-ordering $\preceq^*$ on $\Gamma$ different from $\preceq$ but for which the 
``signs'' of the elements $f_{i,j}$ are the same. To do this, choose $g \!\in\! \Gamma$ such 
that $\Gamma_g$ is different from all of the subgroups $\Gamma_{f_{i,j}}$. This condition implies 
that the corresponding $\Gamma^g$ is different from all of the $\Gamma^{f_{i,j}}$. Now define 
$\preceq'$ as being the extension by $\preceq$ of the extension of the restriction of 
$\preceq$ to $\Gamma_g$ by the restriction of $\bar{\preceq}$ to $\Gamma^g$. One can 
easily show that $\preceq'$ verifies all the desired properties.

Suppose now that, for some integer $n \!\geq\! 1$, there are precisely $n$ subgroups of 
the form $\Gamma_g$. We claim that $\Gamma$ is solvable with solvability length 
at most $n$. Indeed, If $\Gamma_{g_1}$ denotes the maximal proper $\preceq$-convex subgroup 
of $\Gamma$ then, by Conrad's theorem, $\Gamma_{g_1}$ is normal in $\Gamma$, and the 
quotient $\Gamma / \Gamma_{g_1}$ is Abelian. Hence, $\Gamma'$ is contained in $\Gamma_{g_1}$. 
Since $\Gamma_{g_1}$ contains at most $n-1$ subgroups of the form $\Gamma_g$, we may repeat 
this argument... In at most $n$ steps all the $n$-commutators in $\Gamma$ will appear 
to be trivial, which concludes the proof. $\hfill\square$

\vspace{0.4cm}

Left-orderable solvable groups are Conrad-orderable \cite{krop,witte}. Moreover, according to 
\cite{ordering}, if a group $\Gamma$ has infinitely many left-orderings, then no Conradian 
ordering on $\Gamma$ is isolated in $\mathcal{LO}(\Gamma)$. It would be then interesting 
to classify left-orderable solvable groups $\Gamma$ for which $\mathcal{CO}(\Gamma)$ has 
isolated points. 


\section{Bi-orderings on $\mathrm{F}'$}

\hspace{0.35cm} For every dyadic (open, half-open, or closed) interval $I$, we will 
denote by $\mathrm{F}_I$ the subgroup of $\mathrm{F}$ formed by the elements whose support is 
contained in $I$. Notice that if $I$ is closed, then $\mathrm{F}_I$ is isomorphic 
to $\mathrm{F}$. Therefore, 
for every closed dyadic interval $I \! \subset \esp ]0,1[$, every bi-ordering 
$\preceq^*$ on $\mathrm{F}'$ gives 
rise to a bi-ordering on $\mathrm{F} \sim \mathrm{F}_I$. Moreover, if we fix such an $I$, then 
the induced bi-ordering on $\mathrm{F}_I$ completely determines $\preceq^*$ (this is due to the 
invariance by conjugacy). The content of Dlab's theorem consists of the assertion that only a few 
(namely  four) bi-orderings on $\mathrm{F}_I$ may be extended to bi-orderings on $\mathrm{F}'$. 
To reprove this result, we will first focus on a general property of bi-orderings on $\mathrm{F}$.

Let $\preceq$ be a bi-ordering on F. Since bi-invariant orderings are Conradian and $\mathrm{F}$ 
is finitely generated, Conrad's theorem provides us with a (unique up to positive scalar factor) 
non-decreasing group homomorphism $\tau_{\preceq} \!\!: \mathrm{F} \rightarrow (\mathbb{R},+)$ 
whose kernel coincides with the maximal proper $\preceq$-convex subgroup of F. Since 
$\mathrm{F}'$ is a non-Abelian simple group \cite{CFP}, this homomorphism factors 
through $\mathrm{F} / \mathrm{F}' \sim \mathbb{Z}^2$, where the last isomorphism is 
given by \esp 
$f \hspace{0.05cm} \mathrm{F}' \mapsto \big( \log_2 (f'_+(0)), \log_2 (f'_{-}(1)) \big)$. 
\esp Hence, we may write (each representative of the class of) $\tau$ in the form  
$$\tau_{\preceq} (f)= a \log_2 (f'_+ (0)) + b \log_2 (f'_{-}(1)).$$
A canonical representative is obtained by taking $a,b$ so that \esp $a^2 + b^2 = 1$. \esp 
We will call this the {\em normalized Conrad homomorphism} associated to $\preceq$.
In many cases, we will consider this homomorphism as defined on 
$\mathbb{Z}^2 \sim \mathrm{F} / \mathrm{F}'$, so that \esp 
$\tau_{\preceq} \big( (m,n) \big) = am + bn$, and we will 
identify $\tau_{\preceq}$ to the pair $(a,b)$.

Now let $\preceq^*$ be a bi-ordering on $\mathrm{F}'$. For each closed 
dyadic interval $I \! \subset \esp ]0,1[$ let us consider the induced 
bi-ordering on $\mathrm{F} \!\sim\! \mathrm{F}_I$. Since all the subgroups 
$\mathrm{F}_I$ for different closed dyadic intervals are conjugate by 
elements in $\mathrm{F}'$, this induced bi-ordering on $\mathrm{F}$ --which 
we will just denote by $\preceq$-- does not depend on $I$, and hence it 
is inherent to $\preceq^*$. For each such an $I$ let us consider the 
corresponding normalized Conrad homomorphism $\tau_{\preceq,I}$.

\vspace{0.12cm}

\begin{lem} {\em If \esp $\tau_{\preceq}$ corresponds to the pair 
$(a,b)$, then either \esp $a \!=\! 0$ \esp or \esp $b \! = \! 0$.}
\end{lem}

\noindent{\bf Proof.} Assume by contradiction that $a>0$ and $b>0$ (all the other cases 
are analogous). Fix $f \!\in\! \mathrm{F}_{[1/2,3/4]}$ such that $f'_+ (1/2) > 1$ and 
$f'_{-} (3/4) < 1$, and denote $I_1 = [1/4,3/4]$ and $I_2 = [1/2,7/8]$. 
Viewing $f$ as an element in $\mathrm{F}_{I_1} \sim \mathrm{F}$ we have 
$$\tau_{\preceq,I_1} (f) = b \log_2 \big( f'_{-} (3/4)) < 0.$$
Since Conrad's homomorphism is non-decreasing, this implies that $f$ is negative 
with respect to the restriction of $\preceq^*$ to $\mathrm{F}_{I_1}$, and therefore 
$f \prec^* id$. Now viewing $f$ as an element in $\mathrm{F}_{I_2} \sim \mathrm{F}$ we have 
$$\tau_{\preceq,I_2} (f) = a \log_2 \big( f'_{+} (1/2)) > 0,$$
which implies that $f \succ^* id$, thus giving a 
contradiction. $\hfill\square$

\vspace{0.45cm}

We may now pass to the proof of Dlab's theorem. Indeed, assume that for the Conrad's 
homomorphism above one has $a \!>\! 0$ and $b \!=\! 0$. We claim that $\preceq^*$ then 
coincides with $\preceq_{x^{-}}^+$. To show this, we need to show that a non-trivial 
element $f \!\in\! \mathrm{F}'$ is positive with respect to $\preceq^{*}$ if and only 
if $f'_{+} (x_f^{-}) \!>\! 1$. But such an $f$ may be seen as an element in 
$\mathrm{F}_{[x_f^{-},x_f^{+}]}$, and viewed in this way Conrad's homomorphism gives  
$$\tau_{\preceq,[x_f^{-},x_f^{+}]} (f) = a \log_2 (f'_+ (x_f^{-})).$$
Now since $a \!>\! 0$, if $f'_{+} (x_f^{-}) \!>\! 1$ then the right-hand member in this 
equality is positive. Conrad's homomorphism being non-decreasing, this implies that $f$ 
is positive with respect to $\preceq^*$. Analogously, if $f'_{+} (x_f^{-}) \!<\! 1$ 
then $f$ is negative with respect to $\preceq^*$. 

Similar arguments show that the case $a \!<\! 0, \esp b \!=\! 0$ (resp. $a \!=\! 0, \esp b \!>\! 0$, 
and $a \!=\! 0, \esp b \!<\! 0$) necessarily corresponds to the bi-ordering $\preceq_{x^{-}}^{-}$ 
(resp. $\preceq_{x^+}^{-}$, and $\preceq_{x^+}^{+}$), which concludes the proof. 

\vspace{0.1cm}

\begin{qs} According to Proposition \ref{lema-obvio}, a bi-ordering whose positive cone is finitely 
generated as a normal semigroup is completely determined by finitely many inequalities. This 
makes it natural to ask whether this is the case for the restrictions to $\efe'$ of 
$\preceq_{x^{-}}^{+}$, $\preceq_{x^{-}}^{-}$, $\preceq_{x^+}^{+}$, and $\preceq_{x^+}^{-}$. A 
more sophisticated question is the existence of generators \esp $f,g$ \esp of $\efe'$ such that:

\vspace{0.1cm}

\noindent -- $f'_+ (x_f^{-}) > 1$, \esp $g'_+ (x_g^{-}) > 1$, \esp 
$f'_{-} (x_f^{+}) < 1$, \esp and \esp $g'_{-} (x_{g}^{+}) > 1$,

\vspace{0.1cm}

\noindent -- $\efe' \! \setminus \! \{id\}$ is the disjoint union of $\langle \{f,g\} \rangle^+_N$ 
and $\langle \{f^{-1}, g^{-1}\} \rangle^+_{N}$,

\vspace{0.1cm}

\noindent -- $\efe' \! \setminus \! \{id\}$ is also the disjoint union of 
$\langle \{f^{-1},g\} \rangle^+_N$ and $\langle \{f, g^{-1}\} \rangle^+_{N}$.

\vspace{0.1cm}

\noindent A positive answer for the this question would immediately imply Dlab's theorem. 
Indeed, any bi-ordering $\preceq$ on $\efe'$ would be completely determined by the signs  
of $f$ and $g$. For instance, if \esp 
$f \!\succ\! id$ \esp and \esp $g \!\succ\! id$ \esp then $P_{\preceq}^+$ would necessarily 
contain $\langle \{f,g\} \rangle_N^+$, and by the second property above this would 
imply that $\preceq$ coincides with $\preceq_{x^{-}}^{+}$.
\label{prima}
\end{qs}


\section{Bi-orderings on $\mathrm{F}$}

\subsection{Isolated bi-orderings on $\mathrm{F}$}

\hspace{0.35cm} Before classifying all bi-orderings on F, we will first give a proof of the fact 
that the eight elements in $\mathcal{BO}_{Isol} (\mathrm{F})$ are isolated in $\mathcal{BO} 
(\mathrm{F})$. As in the case of $\mathrm{F}'$, this proof strongly uses Conrad's homomorphism.

We just need to consider the cases of $\preceq_{x^{-}}^{+}$ and $\preceq_{0,x^{-}}^{+,-}$. 
Indeed, all the other elements in $\mathcal{BO}_{Isol} (\mathrm{F})$ are obtained from 
these by the action of the (finite Klein's) group generated by the involutions 
$\preceq \hspace{0.05cm} \mapsto \bar{\preceq}$ and 
$\preceq \hspace{0.05cm} \mapsto \preceq_{\sigma}$.

Let us first deal with $\preceq_{x^{-}}^{+}$, denoted $\preceq$ for simplicity. Let 
$(\preceq_k)$ be a sequence in $\mathcal{BO}(\mathrm{F})$ converging to $\preceq$, and let  
$\tau_k \!\sim\! (a_k,b_k)$ be the normalized Conrad's homomorphism for $\preceq_k$ (so 
that \esp $\tau_{k} (m,n) = a_k m + b_k n$ \esp and \esp $a_k^2 + b_k^2 = 1$). \esp 

\vspace{0.35cm}

\noindent{\bf Claim 1.} For $k$ large enough one has $b_k\!=\!0$. 

\vspace{0.1cm}

Indeed, let $f,g$ be two elements in $\mathrm{F}_{]1/2,1]}$ which are positive with respect 
to $\preceq$ and such that $f'_{-}(1)=1/2$ and $g'_{-}(1)=2$. For $k$ large enough, these 
elements must be positive also with respect to $\preceq_k$. Now notice that 
$$\tau_k(f) = -b_k \quad \mbox{ and } \quad \tau_k (g) = b_k.$$
Thus, if $b_k \! \neq \! 0$ then either $f \prec_k id$ or 
$g \prec_k id$, which is a contradiction. Therefore, 
$b_k \! = \! 0$ for $k$ large enough.

\vspace{0.2cm}

Let us now consider the bi-ordering $\preceq^*$ on $\mathrm{F} \sim \mathrm{F}_{[1/2,1]}$ 
obtained as the restriction of $\preceq$. Let $\tau^* \!\sim\! (a^*,b^*)$ 
be the corresponding normalized Conrad's homomorphism. 

\vspace{0.35cm}

\noindent{\bf Claim 2.} One has $b^* \!=\! 0$. 

\vspace{0.1cm}

Indeed, for the elements $f,g$ in $\mathrm{F}_{]1/2,1]}$ above we have 
$$\tau^* (f) = -b^* \quad \mbox{ and } \quad \tau^* (g) = b^*.$$
If $b^* \!\neq\! 0$ this would imply that one of these elements is negative with respect to 
$\preceq^*$, and hence with respect to $\preceq$, which is a contradiction. Thus, $b^* \!=\! 0$.

\vspace{0.2cm}

Denote now by $\preceq_k^*$ the restriction of $\preceq_k$ to $\mathrm{F}_{[1/2,1]}$, and let 
$\tau_k^* \! \sim \! (a_k^*,b_k^*)$ be the corresponding normalized Conrad's homomorphism. 

\vspace{0.35cm}

\noindent{\bf Claim 3.} For $k$ large enough one has $b_k^* \!=\! 0$.

\vspace{0.1cm}

Indeed, the sequence $(\preceq_k^*)$ clearly converges to $\preceq^*$. Knowing 
also that $b^* \!=\! 0$, the proof of this claim is similar to that of Claim 1.

\vspace{0.35cm}

\noindent{\bf Claim 4.} For $k$ large enough one has $a_k \!>\! 0$ and $a_k^* \!>\! 0$.

\vspace{0.2cm}

Since Conrad's homomorphism is non-trivial, both $a_k$ and $a_k^*$ are nonzero. Take any 
$f \!\in\! \mathrm{F}$ such that $f'_+(0) \!=\! 2$. We have $\tau_k (f) \!=\! a_k$. Hence, 
if $a_k < 0$ then \esp $f \prec_k id$, \esp while \esp $f \succ id$... Analogously, if 
$a_k^* < 0$ then one would have \esp $g \prec_k id$ \esp and \esp $g \succ id$ 
\esp for any $g \in \mathrm{F}_{[1/2,1]}$ satisfying \esp $g'(1/2) = 2$.

\vspace{0.35cm}

\noindent{\bf Claim 5.} If $a_k$ and $a_k^*$ are positive and $b_k$ and 
$b_k^*$ are zero, then $\preceq_k$ coincides with $\preceq$.

\vspace{0.2cm}

Given $f \!\in\! \mathrm{F}$ such that $f \succ id$, we need to show that $f$ is positive also 
with respect to $\preceq_k$. If $x_f^{-} \!=\! 0$ then $f'_{+}(0) > 1$, and since $a_k \!>\! 0$ 
this gives $\tau_k (f) = a_k \log_2 (f'_{+}(0)) > 0$, and thus $f \succ_k id$. If 
$x_f^{-} \neq 0$ then $f'_{+}(x_f^{-}) > 1$, and since $a_k^* \!>\! 0$ this gives 
$\tau_k^{*} (f) = a_k^* \log_2 (f'_{+}(x_f^{-})) > 0$, and therefore one still has 
$f \succ_k id$.

\vspace{0.35cm}

The proof for $\preceq_{0,x^{-}}^{+,-}$ is similar to the above one. Indeed, Claims 1, 2, and 3, 
still hold. Concerning Claim 4, one now has that $a_k \!>\! 0$ and $a_k^* \!<\! 0$ for $k$ 
large enough. Having this in mind, one easily concludes that $\preceq_k$ coincides with 
$\preceq_{0,x^{-}}^{+,-}$ for $k$ very large.

\vspace{0.1cm}

\begin{qs} It would be nice to know whether the positive cone of each element in 
$\mathcal{BO}_{Isol} (\mathrm{F})$ is finitely generated as a normal semigroup. 
Notice however that these bi-orderings cannot be completely determined by the 
signs of finitely many elements, since $\mathcal{BO} (\efe)$ is infinite  
(compare Question \ref{prima}).
\label{secua}
\end{qs}


\subsection{Classifying all bi-orderings on $\mathrm{F}$}

\hspace{0.35cm} To simplify, we will denote by $\Lambda$ the union of $\Lambda_{x^{-}}^{+}$, 
$\Lambda_{x^{-}}^{-}$, $\Lambda_{x^{+}}^{+}$, and $\Lambda_{x^{+}}^{-}$. 
To prove our main result, fix a bi-ordering $\preceq$ on F, and let 
$\tau_{\preceq} \!\!: \mathrm{F} \rightarrow (\mathbb{R},+)$ be the corresponding 
normalized Conrad's homomorphism. Since $\tau_{\preceq} \!\sim\! (a,b)$ is 
non-trivial and factors through $\mathbb{Z}^2 \! \sim \! \mathrm{F} / \mathrm{F}'$, 
there are two different cases to be considered. 

\vspace{0.45cm}

\noindent{\bf Case I.} The image $\tau_{\preceq} (\mathbb{Z}^2)$ has rank two.

\vspace{0.15cm}

This case appears when the quotient $a/b$ is irrational. In this case, $\preceq$ induces 
the bi-ordering of irrational type $\preceq_{a/b}$ on $\mathbb{Z}^2$ viewed as 
$\mathrm{F} / \mathrm{F}'$ ({\em c.f.}, Example \ref{abel}). 
Indeed, for each $f \in \mathrm{F} \setminus \mathrm{F}'$ 
the value of $\tau_{\preceq} (f)$ is nonzero, and hence 
it is positive if and only if $f \succ id$. 

The kernel of $\tau_{\preceq}$ coincides with $\mathrm{F}'$. By Dlab's theorem, 
the restriction of $\preceq$ to $\mathrm{F}'$ must coincide with one of the 
bi-orderings $\preceq_{x^{-}}^{+}$, $\preceq_{x^{-}}^{-}$, $\preceq_{x^{+}}^{+}$, 
or $\preceq_{x^{+}}^{-}$. Therefore, $\preceq$ is contained in $\Lambda$, and 
the bi-ordering induced on the $\mathbb{Z}^2$-fiber is of irrational type.

\vspace{0.45cm}

\noindent{\bf Case II.} The image $\tau_{\preceq} (\mathbb{Z}^2)$ has rank one.

\vspace{0.15cm}

This is the difficult case: it appears when either $a/b$ is rational or $b \!=\! 0$. 
There are two sub-cases.

\vspace{0.3cm}

\noindent{\bf Sub-case 1.} Either $a \!=\! 0$ or $b \!=\! 0$.

\vspace{0.15cm}

Assume first that $b \!=\! 0$. Denote by $\preceq^*$ the bi-ordering induced on 
$\mathrm{F}_{[1/2,1]}$, and let $\tau_{\preceq^*} \!\sim\! (a^*,b^*)$ be its 
normalized Conrad's homomorphism. We claim that either $a^{*}$ or $b^*$ is 
equal to zero. Indeed, suppose for instance that $a^* \!>\! 0$ and $b^* \!>\! 0$ 
(all the other cases are analogous). Let $m,n$ be integers such that
\esp $n > 0$ \esp and \esp $a^* m - b^* n > 0,$ \esp 
and let $f$ be an element in $\mathrm{F}_{[3/4,1]}$ 
such that $f'_+(3/4) = 2^m$ and $f'_{-}(1) = 2^{-n}$. 
Then \esp $\tau_{\preceq^*} (f) = -b^* n < 0$, \esp and hence $f \prec id$. 
On the other hand, taking $h \!\in\! \mathrm{F}$ such that $h(3/4) = 1/2$, 
we get that $h^{-1} f h \!\in\! \mathrm{F}_{[1/2,1]}$, and
$$\tau_{\preceq^*} (h^{-1} f h) = a^* \log_2 ((h^{-1}fh)'_+ (1/2)) + 
b^* \log_2 ((h^{-1}fh)'_{-} (1)) = am - bn > 0.$$
But this implies that $h^{-1} f h$, and hence $f$, is positive with respect to 
$\preceq$, which is a contradiction.

\vspace{0.15cm}

\noindent{\bf (i)} If $a \!>\! 0$ and $a^* \! > \! 0$: We claim that $\preceq$ coincides 
with $\preceq_{x^{-}}^+$ in this case. Indeed, let $f \!\in\! \mathrm{F}$ be an element 
which is positive with respect to $\preceq_{x^{-}}^+$. 
We need to show that $f \succ id$. Now, since $a > 0$, if $x_f^{-} \!=\! 0$ then   
$$\tau_{\preceq} (f) = a \log_2 (f'_+ (0)) > 0,$$ 
and hence $f \succ id$. If $x_f^{-} \! \neq \! 0$ then taking $h \!\in\! \mathrm{F}$ such 
that $h (x_f^{-}) \!=\! 1/2$ we obtain that $h^{-1} f h \in \mathrm{F}_{[1/2,1]}$, and
$$\tau_{\preceq^*} (h^{-1} f h) = a^* \log_2 ((h^{-1} f h)'(1/2)) = 
a^* \log_2 (f'(x_f^{-})).$$
Since $a^* \!>\! 0$, the value of the last expression is positive, which implies 
that $h^{-1} f h$, and hence $f$, is positive with respect to $\preceq$.

\vspace{0.15cm}

\noindent{\bf (ii)} If $a \!>\! 0$ and $a^{*} \! < \! 0$: Similar arguments to 
those of (i) above show that $\preceq$ coincides with $\preceq_{0,x^{-}}^{+,-}$ in this case.

\vspace{0.15cm}

\noindent{\bf (iii)} If $a \!>\! 0$ and $b^{*} \!>\! 0$: We claim that $\preceq$ belongs to 
$\Lambda$, and that the induced bi-ordering on the $\mathbb{Z}^2$-fiber is the lexicographic 
one. To show this, we first remark that if $f \in \efe \setminus \efe'$ is positive then either 
$f'_+ (0) \!>\! 1$, or $f'_+ (0) \!=\! 1$ and $f'_{-}(1) \!>\! 1$. Indeed, if $f'_{+}(0) \!\neq\! 1$ 
then the value of \esp $\tau_{\preceq} (f) = a \log_2 (f'_{+}(0)) \neq 0$ \esp must be positive, since 
Conrad's homomorphism is non-decreasing. If $f'_{+} (0) = 1$ we take $h \! \in \! \efe$ such that 
$h (1/2) \!=\! x_f^{-}$. Then $h^{-1}fh$ belongs to $\efe_{[1/2,1]}$, and the value of 
$$\tau_{\preceq^*} (h^{-1}fh) = b^{*} \log_2 ((h^{-1}fh)'_{-}(1)) = b^{*} \log_2 (f'_{-} (1)) \neq 0$$
must be positive, since $f$ (and hence $h^{-1} f h$) is a positive element of $\efe$. 

To show that $\preceq$ induces a bi-ordering on $\mathbb{Z}^2$, we need to check that $\efe'$ is 
$\preceq$-convex. Let $g \!\in\! \efe'$ and $h \!\in\! \efe$ be such that \esp $id \preceq h \preceq g$. 
If $h$ was not contained in $\efe'$, then $hg^{-1}$ would be a negative element in $\efe \setminus \efe'$. 
But since 
$$(hg^{-1})'_{+} (0) = h'_{+}(0) \quad \mbox{ and } \quad 
(hg^{-1})'_{-}(1) = h'_{-}(1),$$
this would contradict the remark above. Therefore, $h$ belongs to $\efe'$, which shows the $\preceq$-convexity 
of $\efe'$. Again, the remark above shows that the induced bi-ordering on $\mathbb{Z}^2$ is the 
lexicographic one.

\vspace{0.15cm}

\noindent{\bf (iv)} If $a \!>\! 0$ and $b^{*} \!<\! 0$: As in (iii) above, $\preceq$ belongs 
to $\Lambda$, and the induced bi-ordering $\preceq_{\mathbb{Z}^2}$ on the $\mathbb{Z}^2$-fiber 
is the one for which $(m,n) \succ_{\mathbb{Z}^2} (0,0)$ if and only if either $m \!>\! 0$, or 
$m \!=\! 0$ and $n \!<\! 0$.

\vspace{0.15cm}

\noindent{\bf (v)} If $a \!<\! 0$ and $a^* \! > \! 0$: As in 
(i) above, $\preceq$ coincides with $\preceq_{0,x^{-}}^{-,+}$ in this case.

\vspace{0.15cm}

\noindent{\bf (vi)} If $a \!<\! 0$ and $a^{*} \! < \! 0$: As 
in (i) above, $\preceq$ coincides with $\preceq_{x^{-}}^{-}$ in this case.

\vspace{0.15cm}

\noindent{\bf (vii)} If $a \!<\! 0$ and $b^{*} \!>\! 0$: As in (iii) above, $\preceq$ belongs 
to $\Lambda$, and the induced bi-ordering $\preceq_{\mathbb{Z}^2}$ on the $\mathbb{Z}^2$-fiber 
is the one for which $(m,n) \succ_{\mathbb{Z}^2} (0,0)$ if and only if either $m \!<\! 0$, or 
$m \!=\! 0$ and $n \!>\! 0$.

\vspace{0.15cm}

\noindent{\bf (viii)} If $a \!<\! 0$ and $b^{*} \!<\! 0$: As in (iii) above, $\preceq$ belongs to 
$\Lambda$, and the induced bi-ordering $\preceq_{\mathbb{Z}^2}$ on the $\mathbb{Z}^2$-fiber is 
the one for which $(m,n) \succ_{\mathbb{Z}^2} (0,0)$ if and only if either $m \!<\! 0$, or 
$m \!=\! 0$ and $n \!<\! 0$.

\vspace{0.15cm}

The case $a \!=\! 0$ is analogous to the preceding one. Letting now $\preceq^*$ 
be the restriction of $\preceq$ to $\mathrm{F}_{[0,1/2]}$, for the normalized Conrad's 
homomorphism $\tau_{\preceq^*} \! \sim \! (a^*,b^*)$ one may check that either 
$a^* \!=\! 0$ or $b^* \! = \! 0$.
 
Assume that $b \!>\! 0$. In the case $b^* \!>\! 0$ (resp. $b^* \!<\! 0$), the 
bi-ordering $\preceq$ coincides with $\preceq_{x^+}^{-}$ (resp. $\preceq_{1,x^+}^{-,+}$). 
If $a^* \!>\! 0$ (resp. $a^* \!<\! 0$), then $\preceq$ corresponds to a point in $\Lambda$ 
whose induced bi-ordering $\preceq_{\mathbb{Z}^2}$ on the $\mathbb{Z}^2$-fiber is the one 
for which \esp $(m,n) \succ_{\mathbb{Z}^2} (0,0)$ \esp if and only if either $n \!>\! 0$, 
or $n \!=\! 0$ and $m \!>\! 0$ (resp. either $n \!>\! 0$, or $n \!=\! 0$ and $m \!<\! 0$).

Assume now that $ b\!<\! 0$. In the case $b^* \!>\! 0$ (resp. $b^* \!<\! 0$), the 
bi-ordering $\preceq$ coincides with $\preceq_{1,x^+}^{+,-}$ (resp. $\preceq_{x^+}^{+}$). 
If $a^* \!>\! 0$ (resp. $a^* \!<\! 0$), then $\preceq$ corresponds to a point in $\Lambda$ 
whose induced bi-ordering $\preceq_{\mathbb{Z}^2}$ on the $\mathbb{Z}^2$-fiber is the one 
for which \esp $(m,n) \succ_{\mathbb{Z}^2} (0,0)$ \esp if and only if either $n \!<\! 0$, 
or $n \!=\! 0$ and $m \!>\! 0$ (resp. either $n \!<\! 0$, or $n \!=\! 0$ and $m \!<\! 0$).

\vspace{0.3cm}

\noindent{\bf Sub-case 2.} Both $a$ and $b$ are nonzero.

\vspace{0.15cm}

The main issue here is to show that $\mathrm{F}'$ is necessarily $\preceq$-convex in $\efe$. But since 
$ker(\tau_{\preceq})$ is already $\preceq$-convex in $\efe$, to prove this it suffices to show 
that $\efe'$ is $\preceq$-convex in $ker (\tau_{\preceq})$. Assume by contradiction that $f$ is 
a positive element in $ker (\tau_{\preceq}) \setminus \efe'$ that is smaller than 
some $h \in \efe'$. Suppose first that $\preceq$ restricted to $\efe'$ coincides with either 
$\preceq_{x^{-}}^{+}$ or $\preceq_{x^{-}}^{-}$, 
and denote by $a$ the leftmost fixed point of $f$ in $]0,1]$. We claim that 
$f$ is smaller than any positive element $g \in \efe_{]0,a[}$. Indeed, since $\preceq$ 
coincides with either $\preceq_{x^{-}}^{+}$ or $\preceq_{x^{-}}^{-}$ 
on $\efe'$, the element $f$ is smaller than any  
positive $\bar{h} \in \efe_{]0,a[}$ such that $x_{\bar{h}}^+$ is to the left of 
$x_h^{-}$; taking $n \in \mathbb{Z}$ such that 
$f^{-n} (x_{h}^{-})$ is to the right of $x_g^{-}$, this gives 
\esp $f = f^{-n} f f^n \prec f^{-n} \bar{h} f^n \prec g.$

Now take a positive element $h_0 \in \efe_{]0,a[}$ such that for \esp $\bar{f}=h_0 f$ 
\esp there is no fixed point in $]0,a[$ (it suffices to consider a positive 
$h_0 \in \efe_{[ \frac{a}{4},\frac{3a}{4}]}$ 
whose graph is very close to the diagonal). Then \esp 
$id \prec \bar{f} \prec h_0 g$ \esp for every positive $g \!\in\! \efe_{]0,a[}$. 
The argument above then shows that $\bar{f}$ is smaller than every positive 
element in $\efe_{]0,a[}$. In particular, since $h_0 = \bar{f} f^{-1}$ is in 
$\efe_{]0,a[}$ and is positive, this implies that \esp $\bar{f} \prec \bar{f} f^{-1}$, 
\esp and hence $f \prec id$, which is a contradiction. 

If the restriction of $\preceq$ to $\efe'$ coincides with either $\preceq_{x^{+}}^+$ 
or $\preceq_{x^{+}}^{-}$, one proceeds similarly but working on the interval $[b,1]$ 
instead of $[0,a]$, where $b$ denotes the rightmost fixed point of $f$ in $[0,1[$. 
This concludes the proof of the $\preceq$-convexity of $\mathrm{F}'$, and hence 
that of our main result.

\vspace{0.1cm}

\begin{rem} Our arguments may be easily modified to show that the subgroup 
$\efe_{-} \!=\! \{f \in \efe \! : f'_{+}(0) = 1 \}$ has six different bi-orderings, 
namely (the restrictions of) $\preceq_{x^{-}}^{+}$, $\preceq_{x^{-}}^{-}$, 
$\preceq_{x^{+}}^{+}$, $\preceq_{x^{+}}^{-}$, $\preceq_{1,x^{+}}^{+,-}$, 
and $\preceq_{1,x^{+}}^{-,+}$. An analogous statement holds for 
$\efe_{+} \!=\! \{f \in \efe \! : f'_{-}(1) = 1 \}$. Finally, the group of 
piecewise-affine orientation-preserving dyadic homeomorphisms of the real line whose 
support is bounded from the right (resp. from the left) admits only two bi-orderings, 
namely (the natural analogues of) $\preceq_{x^{+}}^+$ and $\preceq_{x^{+}}^{-}$ (resp. 
$\preceq_{x^{-}}^+$ and $\preceq_{x^{-}}^{-}$). Notice however that this last 
result is already contained in Dlab's work \cite{dlab} (compare Remark \ref{cuatro}).
\end{rem}


\begin{small}


\vspace{0.13cm}

\noindent Andr\'es Navas\\

\noindent Univ. de Santiago, Alameda 3363, Estaci\'on Central, Santiago, Chile 
(anavas@usach.cl)\\

\vspace{0.2cm}

\noindent Crist\'obal Rivas\\

\noindent Univ. de Chile, Las Palmeras 3425, 
$\tilde{\mathrm{N}}$u$\tilde{\mathrm{n}}$oa, Santiago, Chile 
(cristobalrivas@u.uchile.cl)

\end{small}

\end{document}